\documentclass[12pt]{article}
\usepackage{fullpage}

\usepackage{style}

\DeclareMathOperator{\im}{im}

\def\1{\mathbf{1}}
\DeclareMathOperator\diag{diag}
\DeclareMathOperator\dist{dist}
\DeclareMathOperator{\tr}{trs}

\DeclareMathOperator\IN{in}

\DeclareMathOperator\SNF{SNF}
\DeclareMathOperator\Spec{Spec}
\DeclareMathOperator\SP{sp}

\begin{document}
\title{Distinguishing graphs with two integer matrices}
\author{Carlos A. Alfaro\thanks{
{\tt alfaromontufar@gmail.com}, 
Banco de M\'exico, 
Mexico City, Mexico}
\and
Ralihe R. Villagrán\thanks{{\tt rvillagran@wpi.edu}, Department of Mathematical Sciences, Worcester Polytechnic Institute, Worcester, USA}
\and
Octavio Zapata\thanks{
{\tt octavio@im.unam.mx},
Instituto de Matemáticas,
Universidad Nacional Autónoma de México, 
Mexico City, Mexico}}

\date{}
\maketitle
\begin{abstract}
    It is well known that the spectrum and the Smith normal form of a matrix can be computed in polynomial time.
    Thus, it is interesting to explore how good are these parameters for distinguishing graphs. 
    This is relevant since it is related to the Graph Isomorphism Problem (GIP), which asks to determine whether two graphs are isomorphic.
    In this paper, we explore the computational advantages of using the spectrum and the Smith normal form of two matrices associated with a graph.
    By considering the SNF or the spectrum of two matrices of a graph as a single parameter, we compute the number of non-isomorphic graphs with the same parameter with up to 9 vertices, and with up to 10 vertices when the number of graphs with 9 vertices with the same parameter is less than 1000.
    Focusing on the best 20 combinations of matrices for graphs with 10 vertices, we notice that the number of such graphs with a mate is less than 100 in any of these 20 cases. This computational result improves on similar previous explorations.
    Therefore, the use of the spectrum or the SNF of two matrices at the same time shows a substantial improvement in distinguishing graphs.

\end{abstract}

{\bf Keywords:} walk matrix, degree-distance matrix, transmission-adjacency matrix, Smith normal form, spectrum

{\bf AMS Subject Classification} 05C50

\section{Introduction}

Let $G=(V,E)$ be a simple connected graph with $n$ vertices.
Let $A(G)$ denote the \emph{adjacency matrix} of $G$ and let $\deg(G)$ denote the diagonal matrix with the degrees of the vertices of $G$ in the diagonal.
Given two vertices $u$ and $v$ in $V(G)$, the \emph{distance} $\dist(u,v)$ between $u$ and $v$ is defined as the number of edges in a shortest path joining $u$ and $v$.
Then, the \emph{distance matrix} $D(G)$ of $G$ is the 
$n\times n$ matrix whose $uv$-entry is equal to $\dist(u,v)$.
The \emph{transmission} $\tr(u)$ of vertex $u$ is $\sum_{v\in V(G)}\dist(u,v)$, thus
$\tr(G)$ is the diagonal matrix with the transmissions of the vertices of $G$ in the diagonal.
With the previous definitions, other matrices can be obtained.
The \emph{Laplacian matrix} $L(G)$ is defined as $\deg(G)-A(G)$, the \emph{signless Laplacian matrix} $Q(G)$ is defined as $\deg(G)+A(G)$, the \emph{distance Laplacian matrix} $D^L(G)$ is defined as $\tr(G)-D(G)$, and the \emph{signless distance Laplacian matrix} $D^Q(G)$ is defined as $\tr(G)+D(G)$. 
Furthermore, we define the \emph{degree-distance matrix} $D^{\deg}(G)$ of $G$ as $\deg(G)-D(G)$, the \emph{transmission-adjacency matrix} $A^{\tr}(G)$ of $G$ as $\tr(G)-A(G)$, the \emph{signless degree-distance matrix} $D^{\deg}_+(G)$ of $G$ as $\deg(G)+D(G)$, and the \emph{signless transmission-adjacency matrix} $A^{\tr}_+(G)$ of $G$ as $\tr(G)+A(G)$.
Finally, we define $W_M(G)$ as the matrix $[e\ Me\ \cdots\ M^{n-1}e]$ where $M(G)$ is any of the previous matrices and $e$ is the all-one vector.
In particular, $W_A(G)$ is the \emph{walk matrix}, which is defined as $[e\ A(G)e\ \cdots\ A(G)^{n-1}e]$.


Let us recall the Smith normal form.
Two $n\times n$ matrices $M$ and $N$ are \emph{equivalent} if there exist $P,Q\in GL_n(\mathbb{Z})$ such that $N=PMQ$.
Therefore, if $M$ and $N$ are equivalent, then $M$ can be transformed into $N$ by means of the following operations:
\begin{enumerate}
  \item swap any two rows or any two columns.
  \item add an integer multiple of one row to another row.
  \item add an integer multiple of one column to another column.
  \item multiply any row or column by $\pm 1$.
\end{enumerate}
Thus, if $M$ is a square integer matrix, the Smith normal form (SNF) of $M$ is the unique diagonal matrix $\diag(f_1,f_2,\dots,f_r,0,\dots,0)$ equivalent to $M$, whose non-zero entries are non-negative and satisfy that $f_i$ divides $f_{i+1}$, and $r$ is the rank of $M$. 
The elements $f_1,\dots,f_r$ are known as \emph{invariant factors} of $M$.
The motivation for studying the SNF of matrices, constructed from the combinatorial properties of graphs, is that it gives the description of an Abelian group defined by the cokernel of the matrix.
That is, consider a matrix $M\in\mathbb{Z}^{m\times n}$, the \emph{cokernel} of $M$ is defined as $\mathbb{Z}^{m}/\im(M)$.
Then, the cokernel of $M$ is isomorphic to $\mathbb{Z}_{f_1} \oplus \mathbb{Z}_{f_2} \oplus \cdots \oplus\mathbb{Z}_{f_r} \oplus \mathbb{Z}^{m-r}$,
where $f_1, f_2, \ldots, f_{r}$ are the \emph{invariant factors} of $\rm SNF(M)$ and $r$ is the rank of $M$.
This Abelian group has been investigated under different names, depending on the matrix.
The \emph{Smith group} of $G$ is the cokernel of the adjacency matrix $A(G)$, and the \emph{sandpile group} (or \emph{critical group}) of $G$ is the torsion part of the cokernel of the Laplacian matrix $L(G)$.
For the reader interested in the sandpile group, we recommend the book \cite{Klivans} or the survey \cite{AlfaroMerino}.
Several sandpile groups have been computed, see \cite{alfaval0} and its references.
More recently, the SNF of matrices constructed from the distance matrix has been studied in \cite{BK} and \cite{HW}.
An account of the Smith normal form in combinatorics can be found in \cite{stanley}.

The Graph Isomorphism Problem (GIP) asks to determine whether two graphs are isomorphic, that is to say, whether the adjacency matrices of the two graphs are permutation similar.    
In this sense, if $M=M(G)$ is a matrix constructed from the combinatorial properties of $G$, and then, the eigenvalues of $M$ are analyzed.
The $M$-\emph{spectrum} of $G$ consists of all eigenvalues of $M(G)$. We say that two graphs are $M$-\emph{cospectral} if they have the same $M$-spectrum.
The $M$-spectrum \emph{determines} a graph $G$ if the only graphs which are $M$-cospectral with $G$ are isomorphic to $G$.
Analogously, given two graphs $G$ and $H$, if $M(G)$ and $M(H)$ are equivalent, then, their Smith normal forms are the same.
In such case, $G$ and $H$ are $M$-\emph{coinvariant}. 
The concept of coinvariant graphs was introduced by Vince in \cite{vince}.
Therefore, it is still interesting to find a polynomial-time computable matrix $M$ for which every graph is determined by the spectrum or the SNF of $M$, which solves the graph isomorphism problem. 
In this sense, one of the main questions to tackle is: is it possible to construct a matrix $M$ from the combinatorial properties of the graph such that the $M$-spectrum or its SNF determines all graphs? 
This was backed up by Haemers conjecture: the fraction of graphs with $n$ vertices having an $M$-cospectral mate tends to zero as $n$ tends to infinity.
An affirmative answer was formulated in \cite[Section 2.5]{vDH}. 
Nevertheless, it is more arduous to verify cospectrality of these matrices than to verify isomorphism directly. 

In 1971, Harary, King, Mowshowitz and Read performed a computational enumeration in \cite{HKMR} looking for graphs and digraphs whose adjacency matrices have the same characteristic polynomial, with the purpose of understanding what extent of graphs are determined by the spectrum of the adjacency matrix.
Later, Godsil and McKay \cite{gm} extended the numerical study for cospectral graphs with $n\leq 9$ vertices with respect adjacency matrix. 
Haemers and Spence \cite{HS2004} continued the study for $n = 10, 11$. 
Brouwer and Spence \cite{bs} complemented the study for $n=12$. 
Recently, Aouchiche and Hansen \cite{ah0,ah} also performed computational studies on the cospectrality of matrices obtained from the distance matrix of the connected graphs. 
In \cite{aa}, enumeration results were obtained for the spectrum and the SNF for the matrices $A$, $L$, $Q$, $D$, $D^L$ and $D^Q$.
Recently in \cite{az}, the enumeration results were obtained for the matrices $A^{\tr}$, $A^{\tr}_+$, $D^{\deg}$ and $D^{\deg}_+$. 
In Table~\ref{Tab:cospectralcoinvariantALQDDLDQ}, we recall these findings, which will be described in the following.

\begin{table}[ht]
    \centering
	{\footnotesize
    \begin{tabular}{cccccccc}
	\hline
        $n$ & 4 & 5 & 6 & 7 & 8 & 9 & 10\\
        \hline
        $|\mathcal{G}_n|$ & 6 & 21 & 112 & 853 & 11,117 & 261,080 & 11,716,571\\
        \hline
        $|\mathcal{G}^{sp}_n(A)|$ & 0 & 0 & 2 & 63 & 1,353 & 46,930 & 2,462,141 \\
        $|\mathcal{G}^{sp}_n(L)|$ & 0 & 0 & 4 & 115 & 1,611 & 40,560 & 1,367,215 \\
        $|\mathcal{G}^{sp}_n(Q)|$ & 0 & 2 & 10 & 80 & 1,047 & 17,627 & 615,919  \\
        $|\mathcal{G}^{sp}_n(D)|$ & 0 & 0 & 0 & 22 & 658 & 25,058 & 1,389,986   \\
        $|\mathcal{G}^{sp}_n(D^L)|$ & 0 & 0 & 0 & 43 & 745 & 20,455 & 787,851\\
        $|\mathcal{G}^{sp}_n(D^Q)|$ & 0 & 2 & 6 & 38 & 453 & 8,168 & 319,324\\
        $|\mathcal{G}^{sp}_n(D^{\deg})|$ & 0 & 2 & 6 & 40 & 485 & 9,784 & 355,771 \\
        $|\mathcal{G}^{sp}_n(D^{\deg}_+)|$ & 0 & 0 & 0 & 61 & 901 & 24,095 & 852,504 \\
        $|\mathcal{G}^{sp}_n(A^{\tr})|$ & 0 & 2 & 6 & 38 & 413 & 7,877 & 299,931 \\
        $|\mathcal{G}^{sp}_n(A^{\tr}_+)|$ & 0 & 0 & 0 & 43 & 728 & 19,757 & 765,421 \\
        \hline
        $|\mathcal{G}^{in}_n(A)|$ & 4 & 20 & 112 & 853 & 11,117 & 261,080 & 11,716,571 \\
        $|\mathcal{G}^{in}_n(L)|$ & 2 & 8 & 57 & 526 & 8,027 & 221,834 & 11,036,261 \\
        $|\mathcal{G}^{in}_n(Q)|$ & 2 & 11 & 78 & 620 & 7,962 & 201,282 & 10,086,812 \\
        $|\mathcal{G}^{in}_n(D)|$ & 2 & 15 & 102 & 835 & 11,080 & 260,991 & 11,716,249\\
        $|\mathcal{G}^{in}_n(D^L)|$  & 0 &   0 &   0 &  18 &   455 & 16,505 & 642,002\\
        $|\mathcal{G}^{in}_n(D^Q)|$ & 0 &   2 & 4 & 20 & 259 & 7,444 &  264,955\\
        $|\mathcal{G}^{in}_n(D^{\deg})|$  &  2  &  2  &  6  &  34 & 538  & 17,497  & 902,773 \\
        $|\mathcal{G}^{in}_n(D^{\deg}_+)|$   &  2  & 11  & 46  & 495 & 7,169 & 209,822 & 10,815,879 \\
        $|\mathcal{G}^{in}_n(A^{\tr})|$   &  0  &  2  &  4  & 22  & 240  & 6,642   & 237,118 \\
        $|\mathcal{G}^{in}_n(A^{\tr}_+)|$     &  0  &  0  &  0  & 16  & 456  & 15,952  & 605,625  \\
       
        \hline
	\end{tabular}
	}
\caption{Number of connected graphs with an $M$-cospectral mate and with an $M$-coinvariant mate for $A$, $L$, $Q$, $D$, $D^L$, $D^Q$, $A^{\tr}$, $A^{\tr}_+$, $D^{\deg}$, $D^{\deg}_+$.}
	\label{Tab:cospectralcoinvariantALQDDLDQ}
\end{table}

Let us denote by $\mathcal{G}_n$ the set of connected graphs with $n$ vertices, by $\mathcal{G}^{sp}_n(M)$ the set of graphs in $\mathcal{G}_n$ with a $M$-cospectral mate, and by $\mathcal{G}^{in}_n(M)$ the set of graphs in $\mathcal{G}_n$ with a $M$-coinvariant mate.
The enumeration of the parameters $\mathcal{G}^{sp}_n(M)$ and $\mathcal{G}^{in}_n(M)$ for the matrices $A$, $L$, $Q$, $D$, $D^L$, $D^Q$, $A^{\tr}$, $A^{\tr}_+$, $D^{\deg}$ and $D^{\deg}_+$ is shown in Table~\ref{Tab:cospectralcoinvariantALQDDLDQ}. 
From this data, Abiad and Alfaro \cite{aa} highlighted the potential that the SNF has for distinguishing graphs. 
It seems that this parameter was underestimated since there is a substantial number of connected graphs with a $M$-coinvariant mate for the classic matrices $A$, $L$, $Q$ and $D$, as the reader might note. 
This opens the door to explore more matrices whose SNF might be useful to characterize graphs.
For this reason, in \cite{az} it was explored how good are the spectrum and the SNF of the matrices $A^{\tr}$, $A^{\tr}_+$, $D^{\deg}$ and $D^{\deg}_+$ for distinguish graphs.
The result was that some of these improved previous results.

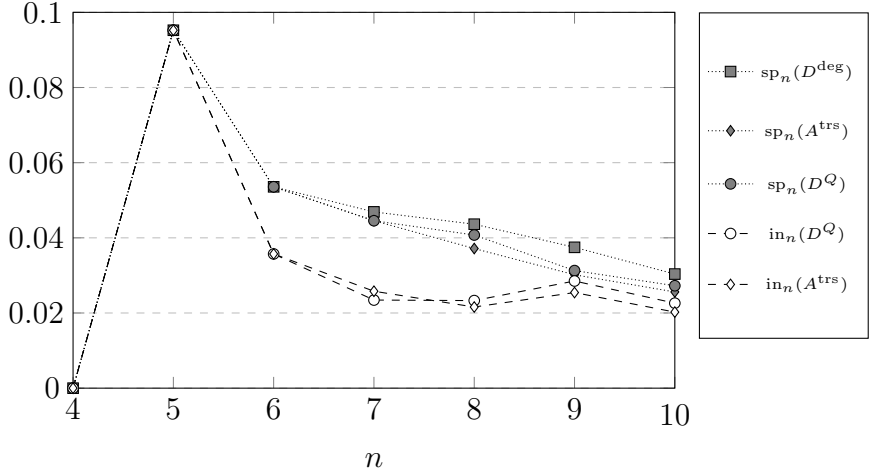
\begin{figure}[ht]
    \centering
    \begin{tikzpicture}[trim axis left]
    \begin{axis}[
        scale only axis,
        title={},
        xlabel={$n$},
        width=8cm, 
        height=5cm,
        xmin=4, xmax=10,
        ymin=0.0, ymax=0.1,
        xtick={4,5,6,7,8,9,10},
        legend style ={ 
            row sep=-0.58cm,
            at={(1.04,1)}, 
            anchor=north west,
            draw=black, 
            fill=white,
            align=left
        },
        ymajorgrids=true,
        grid style=dashed,
        legend columns=1
    ]
    \addplot [
        densely dotted,
        every mark/.append style={solid, fill=gray},
        mark=square*,
        legend entry={\tiny $\SP_n(D^{\deg})$}
        ]
        coordinates {
        (4,0.0) (5,0.09523809523809523) (6,0.05357142857142857) (7,0.04689331770222743) (8,0.04362687775478996) (9,0.037475103416577296) (10,0.030364771399413702)
        };

    \addplot [
        densely dotted,
        every mark/.append style={solid, fill=gray},
        mark=diamond*,
        legend entry={\tiny $\SP_n(A^{\tr})$}
        ]
        coordinates {
        (4,0.0) (5,0.09523809523809523) (6,0.05357142857142857) (7,0.044548651817116064) (8,0.03715031033552217) (9,0.030170828864715796) (10,0.025598871888370754)
        };
        
    \addplot[
        densely dotted,
        every mark/.append style={solid, fill=gray},
        mark=*,
        legend entry={\tiny $\SP_n(D^Q)$}
        ]
        coordinates {
        (4,0.0) (5,0.09523809523809523) (6,0.05357142857142857) (7,0.044548651817116064) (8,0.0407484033462265) (9,0.03128542975333231) (10,0.02725404898754081)
        };
     
     \addplot [
        dashed,
        every mark/.append style={solid, fill=white},
        mark=*,
        legend entry={\tiny $\IN_n(D^Q)$}
        ]
        coordinates {
        (4,0.0) (5,0.09523809523809523) (6,0.03571428571428571) (7,0.023446658851113716) (8,0.023297652244310515) (9,0.02851233338440325) (10,0.022613698154519784)
        };

    \addplot [
        dashed,
        every mark/.append style={solid, fill=white},
        mark=diamond*,
        legend entry={\tiny $\IN_n(A^{\tr})$}
        ]
        coordinates {
        (4,0.0) (5,0.09523809523809523) (6,0.03571428571428571) (7,0.02579132473622509) (8,0.02158855806422596) (9,0.025440478014401715) (10,0.020237832382870382)
        };

    \end{axis}
    \end{tikzpicture}
    \caption{The parameters $sp_n(M)$ and $in_n(M)$ represent the fraction of graphs with $n$ vertices that have at least one $M$-cospectral or $M$-coinvariant mate, respectively. We only show the five best performing parameters.}
    \label{fig:spectruminvariantQDLDQ}
\end{figure}

In order to make it easier to decide which parameter performs better, we introduce the spectral and invariant uncertainty.
The \emph{spectral uncertainty} $\SP_n(M)$  is defined as $|\mathcal{G}^{sp}_n(M)|/|\mathcal{G}_n|$, and similarly, the \emph{invariant uncertainty} $\IN_n(M)$ is defined as $|\mathcal{G}^{in}_n(M)|/|\mathcal{G}_n|$.
In Figure~\ref{fig:spectruminvariantQDLDQ}, these parameters are displayed for the best five parameters in distinguish graphs.
From this, it follows that the SNF of $A^{\tr}$ is better than the SNF and the spectrum of the rest of the matrices. 
This is followed by the SNF of $D^Q$ and by the spectrum of the matrices $A^{\tr}$, $D^Q$ and $D^{\deg}$, respectively.
This was observed in \cite{az} and improves the results obtained in \cite{aa}.

One natural question arises, what happens if we use two parameters simultaneously?
In \cite{aa}, this idea was explored for graphs with 10 vertices and matrices $A$, $L$, $Q$, $D$, $D^L$ and $D^Q$.
In Tables~\ref{Tab:cospectral-invarianttwomatrices} and \ref{Tab:cocospectral-cocoinvarianttwomatrices}, we show results obtained there.
There, we can observe that fewer graphs have a $M$-cospectral and $N$-coinvariant mate at the same time.
For example, only 7,688 graphs with 10 vertices have a $A$-cospectral and $D^Q$-coinvariant mate. 
Which is by far, less than the number of graphs with a cospectral mate or coinvariant mate.
Even, the numbers in the last column of Table~\ref{Tab:cospectralcoinvariantALQDDLDQ} are greater or equal than 200,000.
This suggests that using two parameters at the same time to distinguish graphs could benefit the performance of techniques using these methods for distinguishing graphs.

\begin{table}[h]
	
    \centering
    {
    \footnotesize
    \begin{tabular}{|c|cccccc|}
		\hline
        $M \backslash N$ & $A$ & $L$ & $Q$ & $D$ & $D^L$ & $D^Q$\\
        \hline
        $A$   & 2,151,957 & 24,021 & 22,764 & 1,113,103 & 9,253 & 7,688 \\
        $L$   & 521,200 & 1,059,992 & 121,708 & 192,455 & 562,943 & 44,398 \\
        $Q$   & 136,347 & 84,058 & 486,524 & 48,413 & 44,848 & 250,068 \\
        $D$   & 1,073,185 & 15,176 & 13,496 & 1,145,275 & 8,935 & 7,646 \\
        $D^L$ & 300,596 & 563,219 & 52,757 & 110,574 & 611,989 & 47,004 \\
        $D^Q$ & 65,627 & 44,475 & 245,529 & 28,061 & 46,941 & 255,964 \\
        \hline
    \end{tabular}
    }
	
\caption{$|\mathcal{G}^{sp}_{10}(M)\cap\mathcal{G}^{in}_{10}(N)|$}
	\label{Tab:cospectral-invarianttwomatrices}
\end{table}

\begin{table}[h]
	
    \centering
    {
    \footnotesize
    \begin{tabular}{|c|cccccc|}
		\hline
        $M \backslash N$ & $A$ & $L$ & $Q$ & $D$ & $D^L$ & $D^Q$\\
        \hline
        $A$   &         - & 7,289,194 & 5,712,384 & 11,708,654 & 236,576 &  56,194 \\
        $L$   &    13,864 &         - & 2,150,489 &  4,348,954 & 533,952 &  41,952 \\
        $Q$   &    10,716 &   107,835 &         - &  3,252,128 &  48,771 & 237,195 \\
        $D$   & 1,188,225 &     9,837 &     7,762 &          - &  92,915 &  24,517 \\
        $D^L$ &     9,823 &   722,060 &    58,557 &      9,449 &       - &  46 239 \\
        $D^Q$ &     7,726 &    57,897 &   309,212 &      7,712 &  61,909 &       - \\
        \hline
    \end{tabular}
    }
	
\caption{The lower triangular entries shows $|\mathcal{G}^{sp}_{10}(M)\cap\mathcal{G}^{sp}_{10}(N)|$ and the upper triangular entries shows $|\mathcal{G}^{in}_{10}(M)\cap\mathcal{G}^{in}_{10}(N)|$}
	\label{Tab:cocospectral-cocoinvarianttwomatrices}
\end{table}

In this paper, we continue with this exploration.
First, we extend the previous computational enumeration to all connected graphs with at most 10 vertices with a cospectral mate with respect to the matrices $W_A$, $W_D$, $W_Q$, $W_{D^Q}$, $W_{A^{\tr}}$ and $W_{D^{\deg}}$.
We will see that the spectrum of these matrices has much better performance than the matrices $A$, $L$, $Q$, $D$, $D^L$, $D^Q$, $A^{\tr}$, $A^{\tr}_+$, $D^{\deg}$ and $D^{\deg}_+$.
This is not the case for the spectrum and SNF of the rest of the $W_M$-type matrices, whose invariant uncertainty for connected graphs with 9 vertices is above 0.224. 
Another point of interest is whether the spectrum or the SNF of the $W_M$ matrices can used to distinguish trees.
Surprisingly, only for the $W_A$-matrix there are no $W_A$-cospectral trees up to 20 vertices.
The other parameters of the $W_M$-type matrices are not good for distinguishing trees.
Finally, we turn our attention to the use of two parameters at the same time to distinguish graphs.
Enumeration results show an impressive performance, for instance, 20 combinations of two parameters show less than 100 connected graphs with 10 vertices having a mate, see Table~\ref{tab:best20}.
In particular, there are only 80 connected graphs with 10 vertices having a $W_A$-cospectral and $D^L$-coinvariant mate at the same time.

\section{Experiments}

\subsection{$W_M$-matrices}

Firstly, we explore the $W_M$-type matrices. First note that the matrix $W_L$ ($W_{D^L}$) is not of interest since the vector with all entries equal to 1 is in the kernel of the (distance) Laplacian matrix.
The spectral and invariant uncertainty of the $W_M$ matrices of graphs with 9 vertices is shown in Table~\ref{tab:Wtypematrices}.

\begin{table}[h]
    \centering
    {
    \footnotesize
    \begin{tabular}{cccccccccc}
        \hline
         $M$ & $A$ & $D$ & $Q$  & $D^Q$ & $A^{\tr}$ & $A^{\tr}_+$ & $D^{\deg}$ & $D^{\deg}_+$ \\
        \hline
        $sp_9(W_M)$ & 0.0071  & 0.0048 & 0.0071 & 0.0049 & 0.0048 & 0.4096 & 0.0048 & 0.4097  \\
        $in_9(W_M)$ & 0.9987 & 0.5827 & 0.9712 & 0.2662 & 0.2240 & 0.7198 & 0.2303 & 0.7245 \\
        \hline
    \end{tabular}
    }
    \caption{The 9-\emph{th} spectral and invariant uncertainty of $W_M$-type matrices.}
    \label{tab:Wtypematrices}
\end{table}

From Table~\ref{tab:Wtypematrices}, it follows that only the spectrum of the matrices $W_Q$, $W_D$, $W_{D^Q}$, $W_{A^{\tr}}$ have possibility of distinguish graphs better than previous results.
They are followed, in performance, by the SNF of $W_{D^Q}$ and $W_{A^{\tr}}$, whose uncertainty is not smaller than $0.1$, for this reason, we will omit them in Table~\ref{tab:spectrumWtypematrices}.
Note that since $W_L$ and $W_{D^L}$ are equal to the zero matrix with the first column of ones, then $sp_9(W_L)$, $sp_9(W_{D^L})$, $in_9(W_L)$ and $in_9(W_{D^L})$ are equal to 1.
Table~\ref{tab:spectrumWtypematrices} shows the spectral uncertainty of the matrices $W_Q$, $W_D$, $W_{D^Q}$, $W_{A^{\tr}}$ up to $n=10$.

\begin{table}[h]
    \centering
    {
    \footnotesize
    \begin{tabular}{cccccccc}
        \hline
        $n$ & 4 & 5 & 6 & 7 & 8 & 9 & 10\\
        \hline
        $|\mathcal{G}^{sp}_n(W_A)|$ & 0 & 0 & 6 & 20 & 191 & 1,868 & 39,105 \\
        $|\mathcal{G}^{sp}_n(W_D)|$ & 0 & 0 & 4 & 16 & 124 & 1,279 & 28,784 \\
        $|\mathcal{G}^{sp}_n(W_Q)|$ & 0 & 0 & 6 & 20 & 191 & 1,874 & 39,286 \\
        $|\mathcal{G}^{sp}_n(W_{D^Q})|$ & 0 & 0 & 4 & 16 & 124 & 1,280 & 28,789 \\
        $|\mathcal{G}^{sp}_n(W_{A^{\tr}})|$ & 0 & 0 & 4 & 16 & 120 & 1,278 & 28,651 \\
        $|\mathcal{G}^{sp}_n(W_{D^{\deg}})|$ & 0 & 0 & 4 & 16 & 120 & 1,278 & 28,669 \\
        \hline
    \end{tabular}
    }
    \caption{Number of connected graphs with an $M$-cospectral mate for $W_A$, $W_Q$, $W_D$, $W_{D^Q}$, $W_{A^{\tr}}$ and $W_{D^{\deg}}$.}
    \label{tab:spectrumWtypematrices}
\end{table}

\begin{figure}[ht]
    \centering
    \begin{tikzpicture}[trim axis left]
    \begin{axis}[
        scale only axis,
        title={},
        xlabel={$n$},
        width=8cm, 
        height=8cm,
        xmin=6, xmax=10,
        ymin=0.0, ymax=0.055,
        xtick={4,5,6,7,8,9,10},
        legend style ={ 
            row sep=-0.71cm,
            at={(1.04,1)}, 
            anchor=north west,
            draw=black, 
            fill=white,
            align=left
        },
        ymajorgrids=true,
        grid style=dashed,
        legend columns=1
    ]
    \addplot [
        densely dotted,
        every mark/.append style={solid, fill=gray},
        mark=square*,
        legend entry={\tiny $\SP_n(D^{\deg})$}
        ]
        coordinates {
        (4,0.0) (5,0.09523809523809523) (6,0.05357142857142857) (7,0.04689331770222743) (8,0.04362687775478996) (9,0.037475103416577296) (10,0.030364771399413702)
        };

    \addplot [
        densely dotted,
        every mark/.append style={solid, fill=gray},
        mark=diamond*,
        legend entry={\tiny $\SP_n(A^{\tr})$}
        ]
        coordinates {
        (4,0.0) (5,0.09523809523809523) (6,0.05357142857142857) (7,0.044548651817116064) (8,0.03715031033552217) (9,0.030170828864715796) (10,0.025598871888370754)
        };
        
    \addplot[
        densely dotted,
        every mark/.append style={solid, fill=gray},
        mark=*,
        legend entry={\tiny $\SP_n(D^Q)$}
        ]
        coordinates {
        (4,0.0) (5,0.09523809523809523) (6,0.05357142857142857) (7,0.044548651817116064) (8,0.0407484033462265) (9,0.03128542975333231) (10,0.02725404898754081)
        };
     
     \addplot [
        dashed,
        every mark/.append style={solid, fill=white},
        mark=*,
        legend entry={\tiny $\IN_n(D^Q)$}
        ]
        coordinates {
        (4,0.0) (5,0.09523809523809523) (6,0.03571428571428571) (7,0.023446658851113716) (8,0.023297652244310515) (9,0.02851233338440325) (10,0.022613698154519784)
        };

    \addplot [
        dashed,
        every mark/.append style={solid, fill=white},
        mark=diamond*,
        legend entry={\tiny $\IN_n(A^{\tr})$}
        ]
        coordinates {
        (4,0.0) (5,0.09523809523809523) (6,0.03571428571428571) (7,0.02579132473622509) (8,0.02158855806422596) (9,0.025440478014401715) (10,0.020237832382870382)
        };

    \addplot [
        densely dotted,
        every mark/.append style={solid, fill=gray},
        mark=x,
        legend entry={\tiny $\SP_n(W_Q)$}
        ]
        coordinates {
        (4,0)	(5,0)	(6,0.0535714285714286)	(7,0.0234466588511137)	(8,0.0171808941261132)	(9,0.00717787651294622)	(10,0.00335302879997911)
        };

    \addplot [
        densely dotted,
        every mark/.append style={solid, fill=gray},
        mark=+,
        legend entry={\tiny $\SP_n(W_A)$}
        ]
        coordinates {
        (4,0)	(5,0)	(6,0.0535714285714286)	(7,0.0234466588511137)	(8,0.0171808941261132)	(9,0.00715489505132526)	(10,0.00333758059418579)
        };

    \addplot [
        densely dotted,
        every mark/.append style={solid, fill=gray},
        mark=Mercedes star,
        legend entry={\tiny $\SP_n(W_{D^Q})$}
        ]
        coordinates {
        (4,0)	(5,0)	(6,0.0357142857142857)	(7,0.018757327080891)	(8,0.0111540883331834)	(9,0.00490271181247127)	(10,0.00245711821317005)
        };

    \addplot [
        densely dotted,
        every mark/.append style={solid, fill=gray},
        mark=|,
        legend entry={\tiny $\SP_n(W_D)$}
        ]
        coordinates {
        (4,0)	(5,0)	(6,0.0357142857142857)	(7,0.018757327080891)	(8,0.0111540883331834)	(9,0.00489888156886778)	(10,0.00245669146715366)
        };

    \addplot [
        densely dotted,
        every mark/.append style={solid, fill=gray},
        mark=Mercedes star flipped,
        legend entry={\tiny $\SP_n(W_{D^{\deg}})$}
        ]
        coordinates {
        (4,0)	(5,0)	(6,0.0357142857142857)	(7,0.018757327080891)	(8,0.010794279032113)	(9,0.00489505132526429)	(10,0.00244687630877669)
        };

    \addplot [
        densely dotted,
        every mark/.append style={solid, fill=gray},
        mark=halfsquare*,
        legend entry={\tiny $\SP_n(W_{A^{\tr}})$}
        ]
        coordinates {
        (4,0)	(5,0)	(6,0.0357142857142857)	(7,0.018757327080891)	(8,0.010794279032113)	(9,0.00489505132526429)	(10,0.00244534002311769)
        };
        
    \end{axis}
    \end{tikzpicture}
    \caption{The eleven best matrices whose spectral or invariant uncertainties have lower values for connected graphs with up to 10 vertices.}
    \label{fig:spectruminvariantWs}
\end{figure}
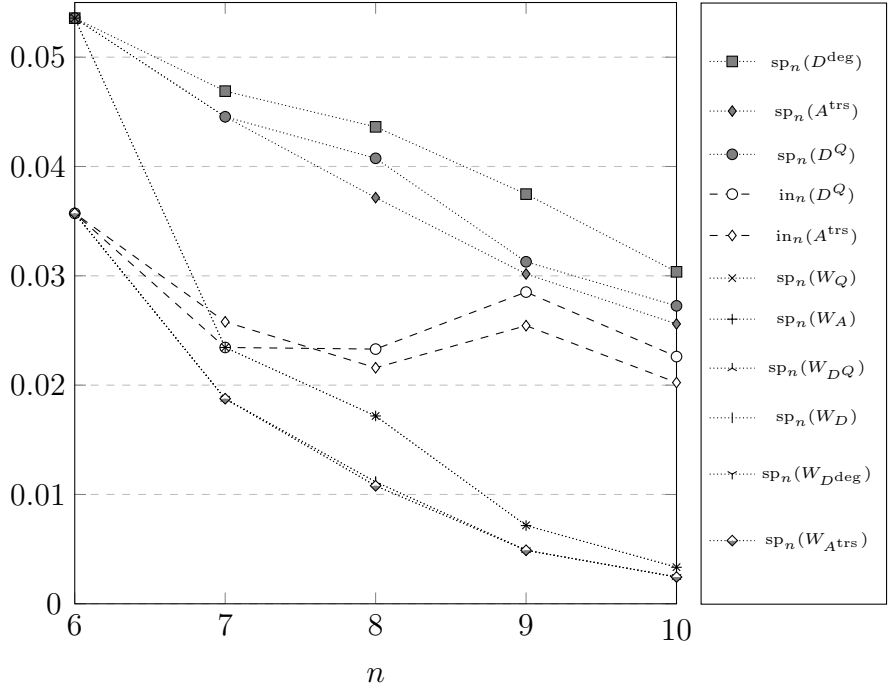

In Figure~\ref{fig:spectruminvariantWs}, we use this information to update Figure~\ref{fig:spectruminvariantQDLDQ}.
There is a clear advantage in using the spectrum of the matrices $W_Q$, $W_D$, $W_{D^Q}$, $W_{A^{\tr}}$ to distinguish graphs comparing to the best parameters previously studied.

\subsection{$W_M$-cospectral trees}

Particularly on trees, we recall the following.
In \cite{ah0}, Aouchiche and Hansen reported on $D$, $D^L$ and $D^Q$ cospectral trees with up to 20 vertices.
For matrix $D$, there were found only 22 pairs of $D$-cospectral mates in all 1,346,023 trees with up to 20 vertices.
Among the 317,955 trees on 19 vertices, there are six pairs of $D$-cospectral mates.
There are 14 pairs of $D$-cospectral mates over all the 823,065 trees on 20 vertices.
Outstandingly, there are no $D^L$-cospectral mates or $D^Q$-cospectral mates amongst these trees.
Considering this, Aouchiche and Hansen conjectured that the spectrum of $D^L$ and $D^Q$ determines every tree.
In the same direction, Abiad and Alfaro \cite{aa} reported enumeration results on $M$-coinvariant trees for the same matrices.
For matrix $D$, they proved that all trees with $n$ vertices are $D$-coinvariant.
This result follows from a result of Hou and Woo \cite{HW} that obtains the Smith normal form of the distance matrix of a tree, which extends Graham and Pollak celebrated formula on the determinant of the distance matrix of a tree.
On the other hand, Abiad and Alfaro reported that there are no $D^L$-coinvariant mates or $D^Q$-coinvariant mates among all trees with up to 20 vertices.
Thus, Abiad and Alfaro conjectured that the SNF of the matrices $D^L$ and $D^Q$ determines every tree.
Recently, in \cite{az}, the authors extended previous enumeration results to the matrices $A^{\tr}$, $A^{\tr}_+$, $D^{\deg}$ and $D^{\deg}_+$ of trees with up to 20 vertices.
The findings were the following.
There are no $D^{\deg}_+$-cospectral trees nor $D^{\deg}$-cospectral trees nor $A^{\tr}$-cospectral trees with up to 20 vertices.
However, there are four $A^{\tr}_+$-cospectral trees with 20 vertices and no $A^{\tr}_+$-cospectral trees with up to 19 vertices.
Considering the SNF, there are no $A^{\tr}_+$-coinvariant trees nor $A^{\tr}$-coinvariant trees with up to 20 vertices.
However, we found two $D^{\deg}$-coinvariant trees with 14 vertices,
6 $D^{\deg}$-coinvariant trees with 16 vertices,
14 $D^{\deg}$-coinvariant trees with 17 vertices,
22 $D^{\deg}$-coinvariant trees with 18 vertices,
40 $D^{\deg}$-coinvariant trees with 19 vertices, and 
76 $D^{\deg}$-coinvariant trees with 20 vertices.
With respect to the matrix $D^{\deg}_+$, we found two $D^{\deg}_+$-coinvariant trees with 9 vertices, 
6 $D^{\deg}_+$-coinvariant trees with 10 vertices, 
20 $D^{\deg}_+$-coinvariant trees with 11 vertices, 
46 $D^{\deg}_+$-coinvariant trees with 12 vertices, 
148 $D^{\deg}_+$-coinvariant trees with 13 vertices, 
373 $D^{\deg}_+$-coinvariant trees with 14 vertices, 
1093 $D^{\deg}_+$-coinvariant trees with 15 vertices, 
2912 $D^{\deg}_+$-coinvariant trees with 16 vertices, 
8189 $D^{\deg}_+$-coinvariant trees with 17 vertices, 
22551 $D^{\deg}_+$-coinvariant trees with 18 vertices, 
64738 $D^{\deg}_+$-coinvariant trees with 19 vertices, and 
183211 $D^{\deg}_+$-coinvariant trees with 20 vertices.
Considering this evidence, the authors conjectured in \cite{az} that trees are determined by spectrum of the matrices $D^{\deg}_+$, $D^{\deg}$ and $A^{\tr}$, and also determined by the SNF of the matrices $A^{\tr}_+$ and $A^{\tr}$.

Now in the context of the $W_M$-type matrices, we explore cospectral trees for $W_A$, $W_D$, $W_Q$, $W_{D^Q}$, $W_{A^{\tr}}$ and $W_{D^{\deg}}$.
There are no $W_A$-cospectral trees up to 20 vertices.
For matrix $W_D$, we found no $W_D$-cospectral trees with up to 18 vertices, but 3 and 7 $W_D$-cospectral trees with $19$ and $20$ vertices, respectively.
For matrix $W_Q$, we found no $W_Q$-cospectral trees with up to 19 vertices and 7 $W_Q$-cospectral trees with 20 vertices.
For matrix $W_{D^Q}$, we found no $W_{D^Q}$-cospectral trees with up to 18 vertices, and 2 and 7 $W_{D^Q}$-cospectral trees with 19 and 20 vertices, respectively.
For matrix $W_{A^{\tr}}$, we found no $W_{A^{\tr}}$-cospectral trees with up to 18 vertices, ad 4 and 15 $W_{A^{\tr}}$-cospectral trees with 19 and 20 vertices, respectively.
For matrix $W_{D^{\deg}}$, we found no $W_{D^{\deg}}$-cospectral trees up to 19 vertices, and 11 $W_{D^{\deg}}$-cospectral trees with 20 vertices.

\subsection{Distinguishing graphs with two parameters}

Now, we turn our attention to distinguishing graphs with the SNF or the spectrum of two matrices.
By considering the SNF or the spectrum of two matrices as a unique parameter, we compute the number of non-isomorphic graphs with this same parameter with up to 9 vertices, and with up to 10 vertices when the number of graphs with 9 vertices with the same parameter is less than 1000.
The complete results are shown in the tables of the Appendix and at \url{https://github.com/alfaromontufar/DistinguishingGraphs} in {\tt xlsx} format.
We focused on these graphs since they show the best behavior, and obtaining this information requires a lot of computational consumption. 
Therefore, from Table~\ref{tab:datos5}, some values for $n=10$ are unknown, which are indicated by a dash mark.
The software we used to make the above computations was SageMath 9.5 \cite{sage} on a MacBook Pro (16-inch, 16GB RAM, 1TB Storage, 2.3GHz Intel Core i9) machine.
Considering this information, in Table~\ref{tab:best20} we highlight the 20 best combinations of spectrum and SNF of matrices for graphs with 10 vertices.
It is impressive to note that in such cases the number of graphs with 10 vertices having a mate is less than 100!

\begin{table}[h]
    \centering
    {
    \footnotesize
    \begin{tabular}{c@{\extracolsep{1cm}}c}
        \begin{tabular}{cc}
            \hline
            parameter & \# graphs\\
            \hline
            $\Spec(W_A)\cap\SNF(D^L)$ &	80\\
            $\Spec(W_D)\cap\SNF(D^L)$ &	80\\
            $\Spec(A^{\tr})\cap\SNF(D^L)$ &	82\\
            $\Spec(D^{\deg})\cap\SNF(D^L)$ &	82\\
            $\Spec(W_A)\cap\SNF(A^{\tr}_+)$ &	82\\
            $\Spec(W_D)\cap\SNF(A^{\tr}_+)$ &	82\\
            $\Spec(A^{\tr})\cap\SNF(A^{\tr}_+)$ &	84\\
            $\Spec(D^{\deg})\cap\SNF(A^{\tr}_+)$ &	84\\
            $\Spec(W_{D^Q})\cap\SNF(A^{\tr}_+)$ &	84\\
            $\Spec(W_{D^Q})\cap\SNF(D^L)$ &	84\\
            \hline
        \end{tabular}
        &
        \begin{tabular}{cc}
            \hline
            parameter & \# graphs\\
            \hline
            $\Spec(W_Q)\cap\SNF(A^{\tr}_+)$ &	84\\
            $\Spec(W_Q)\cap\SNF(D^L)$ &	84\\
            $\Spec(A^{\tr})\cap\SNF(L)$ &	86\\
            $\Spec(D^{\deg})\cap\SNF(L)$ &	86\\
            $\Spec(W_D)\cap\SNF(L)$ &	88\\
            $\Spec(W_{D^Q})\cap\SNF(L)$ &	92\\
            $\Spec(W_A)\cap\Spec(D^L)$ &	97\\
            $\Spec(A^{\tr})\cap\Spec(D^L)$ &	99\\
            $\Spec(D^{\deg})\cap\Spec(D^L)$ &	99\\
            $\Spec(W_D)\cap\Spec(D^L)$ &	99\\
            \hline
        \end{tabular}
    \end{tabular}
    }
    \caption{Number of graphs with $10$ vertices having a mate with respect to the spectrum or the SNF of two matrices associated with the graph at the same time.}
    \label{tab:best20}
\end{table}

The best results were obtained by combining the SNF of $D^L$ and the spectrum of $W_A$ and $W_D$.
Only 80 connected graphs with 10 vertices have a mate with such a parameter.
This strengthens the use of the SNF together with the spectrum of the $W_M$-type matrices.
The next best results obtained only 82 graphs with a mate.
It is interesting to note that in these cases the $A^{\tr}$, $D^{\deg}$ or the $A^{\tr}_+$ matrices appeared.
Another point of interest is that in positions 13 to 16, the SNF of the Laplacian matrix appeared, this is interesting since we already have observed that this parameter alone is not good for distinguishing graphs, as can be observed in Table~\ref{Tab:cospectralcoinvariantALQDDLDQ}, and that the SNF of the Laplacian matrix is used to compute the algebraic structure of the sandpile group.

\begin{figure}
    \centering
    \begin{tabular}{c@{\extracolsep{1cm}}c}
    \begin{tikzpicture}[scale=1,thick]
        \tikzstyle{every node}=[minimum width=0pt, inner sep=2pt, circle]
        \draw (0:1) node[draw] (6) {\tiny };
        \draw (36:1) node[draw] (0) {\tiny };
        \draw (72:1) node[draw] (3) {\tiny };
        \draw (108:1) node[draw] (8) {\tiny };
        \draw (144:1) node[draw] (2) {\tiny };
        \draw (180:1) node[draw] (5) {\tiny };
        \draw (216:1) node[draw] (9) {\tiny };
        \draw (252:1) node[draw] (4) {\tiny };
        \draw (288:1) node[draw] (7) {\tiny };
        \draw (324:1) node[draw] (1) {\tiny };
        \draw  (0) edge (3);
        \draw  (0) edge (4);
        \draw  (0) edge (6);
        \draw  (0) edge (7);
        \draw  (0) edge (8);
        \draw  (1) edge (4);
        \draw  (1) edge (5);
        \draw  (1) edge (6);
        \draw  (1) edge (7);
        \draw  (1) edge (9);
        \draw  (2) edge (5);
        \draw  (2) edge (6);
        \draw  (2) edge (7);
        \draw  (2) edge (8);
        \draw  (2) edge (9);
        \draw  (3) edge (5);
        \draw  (3) edge (6);
        \draw  (3) edge (7);
        \draw  (3) edge (8);
        \draw  (4) edge (7);
        \draw  (4) edge (8);
        \draw  (4) edge (9);
        \draw  (5) edge (8);
        \draw  (5) edge (9);
        \draw  (6) edge (9);
        \draw (0,-1.5) node () { \text{\texttt{ICpvfq{\char`\{}Z{\char`\_}}} };
    \end{tikzpicture}
    &
    \begin{tikzpicture}[scale=1,thick]
        \tikzstyle{every node}=[minimum width=0pt, inner sep=2pt, circle]
        \draw (0:1) node[draw] (0) {\tiny };
        \draw (36:1) node[draw] (4) {\tiny };
        \draw (72:1) node[draw] (9) {\tiny };
        \draw (108:1) node[draw] (2) {\tiny };
        \draw (144:1) node[draw] (5) {\tiny };
        \draw (180:1) node[draw] (3) {\tiny };
        \draw (216:1) node[draw] (7) {\tiny };
        \draw (252:1) node[draw] (1) {\tiny };
        \draw (288:1) node[draw] (8) {\tiny };
        \draw (324:1) node[draw] (6) {\tiny };
        \draw  (0) edge (3);
        \draw  (0) edge (4);
        \draw  (0) edge (6);
        \draw  (0) edge (7);
        \draw  (0) edge (8);
        \draw  (1) edge (4);
        \draw  (1) edge (5);
        \draw  (1) edge (6);
        \draw  (1) edge (7);
        \draw  (1) edge (8);
        \draw  (2) edge (4);
        \draw  (2) edge (5);
        \draw  (2) edge (6);
        \draw  (2) edge (7);
        \draw  (2) edge (9);
        \draw  (3) edge (5);
        \draw  (3) edge (7);
        \draw  (3) edge (8);
        \draw  (3) edge (9);
        \draw  (4) edge (8);
        \draw  (4) edge (9);
        \draw  (5) edge (7);
        \draw  (5) edge (9);
        \draw  (6) edge (8);
        \draw  (6) edge (9);
        \draw (0,-1.5) node () { \text{\texttt{ICxvFjYN{\char`\_}}} };
    \end{tikzpicture}
    \end{tabular}
    \caption{A pair of $\Spec(W_A)\cap\SNF(D^L)$-mates with 10 vertices.}
    \label{fig:cospectralcoinvariantmates}
\end{figure}
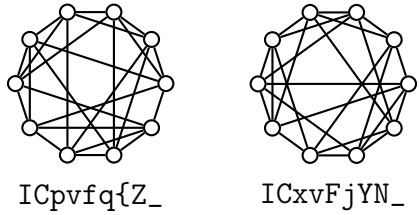

The list of mates of such 20 cases can be found in separate files at \url{https://github.com/alfaromontufar/DistinguishingGraphs}.
Each row in the files corresponds to a set of mates in {\tt graph6} format.
For instance, the file \texttt{CoSpec-WA\_and\_CoSNF-DL\_10.txt} contains in each row the $\Spec(W_A)\cap\SNF(D^L)$-mates with 10 vertices.
And, the graphs in the first row are the $\Spec(W_A)\cap\SNF(D^L)$-mates: \texttt{ICpvfq\{Z\_} and \texttt{ICxvFjYN\_}, which are shown in Figure \ref{fig:cospectralcoinvariantmates}.

\section{Conclusions}

We are interested in the use of the spectrum or SNF of two matrices associated with a graph at the same time as a method to distinguish graphs.
This is relevant since it is related to the famous Graph Isomorphism Problem (GIP), which asks to determine whether two graphs are isomorphic.
Since computing the spectrum or the SNF can be achieved in polynomial time, if some of these invariants distinguish all graphs, then the GIP is polynomial time.

In particular, we introduce $W_M(G)$ as the matrix $[e\ Me\ \cdots\ M^{n-1}e]$ where $M(G)$ can be any matrix associated with $G$, and $e$ is the all-one vector.
In particular, $W_A(G)$ is the walk matrix.
We found that the spectrum of the matrices $W_Q$, $W_D$, $W_{D^Q}$, $W_{A^{\tr}}$ and $W_{D^{\deg}}$ are good for distinguishing connected graphs with up to 10 vertices.
In particular, they perform better than previously studied parameters, which are shown in Table~\ref{Tab:cospectralcoinvariantALQDDLDQ}.
This is summarized in Figure~\ref{fig:spectruminvariantWs}. 
However, only for the $W_A$ matrix, there are no $W_A$-cospectral trees up to 20 vertices.
The other $W_M$ matrices have at least a pair of cospectral trees with less or equal to 20 vertices.

Finally, by considering the SNF or the spectrum of two matrices of a graph as a unique parameter, we compute the number of non-isomorphic graphs with this same parameter with up to 9 vertices, and with up to 10 vertices when the number of graphs with 9 vertices with the same parameter is less than 1000.
Focusing on the 20 best combinations of matrices for graphs with 10 vertices, it is impressive to note that in such cases the number of graphs with 10 vertices having a mate is less than 100.
Therefore, the use of the spectrum or the SNF of two matrices at the same time shows a substantial improvement in distinguishing graphs.
Which will propel new theoretical results in distinguishing graphs using two parameters at the same time.

\section*{Acknowledgement}
The research of C.A. Alfaro is partially supported by the Sistema Nacional de Investigadores grant number 220797. 
The research of O. Zapata is partially supported by the Sistema Nacional de Investigadores grant number 620178.

\bibliographystyle{plain}
\bibliography{bibliography}

\appendix
\section{Tables}

\begin{table}[h]
    \footnotesize
    \centering
\begin{tabular}{cccccc}
parameter & 6 & 7 & 8 & 9 & 10 \\
\hline
$|\mathcal{G}^{sp}_n(W_A)\cap\mathcal{G}^{in}_n(D^L)|$ & 0 & 0 & 0 & 0 & 80 \\
$|\mathcal{G}^{sp}_n(W_D)\cap\mathcal{G}^{in}_n(D^L)|$ & 0 & 0 & 0 & 0 & 80 \\
$|\mathcal{G}^{sp}_n(W_A)\cap\mathcal{G}^{in}_n(A^{\tr}_+)|$ & 0 & 0 & 0 & 0 & 82 \\
$|\mathcal{G}^{sp}_n(W_D)\cap\mathcal{G}^{in}_n(A^{\tr}_+)|$ & 0 & 0 & 0 & 0 & 82 \\
$|\mathcal{G}^{sp}_n(W_{A^{\tr}})\cap\mathcal{G}^{in}_n(D^L)|$ & 0 & 0 & 0 & 0 & 82 \\
$|\mathcal{G}^{sp}_n(W_{D^{\deg}})\cap\mathcal{G}^{in}_n(D^L)|$ & 0 & 0 & 0 & 0 & 82 \\
$|\mathcal{G}^{sp}_n(W_Q)\cap\mathcal{G}^{in}_n(D^L)|$ & 0 & 0 & 0 & 0 & 84 \\
\end{tabular}
    \caption{Number of graphs with up to $10$ vertices having a mate with respect to the spectrum or the SNF of two matrices associated with the graph at the same time. (part 1)}
    \label{tab:datos1}
\end{table}

\begin{table}[h]
    \footnotesize
    \centering
    \begin{tabular}{cccccc}
parameter & 6 & 7 & 8 & 9 & 10 \\
\hline
$|\mathcal{G}^{sp}_n(W_D)\cap\mathcal{G}^{in}_n(A^{\tr})|$ & 0 & 0 & 0 & 2 & 284 \\
$|\mathcal{G}^{sp}_n(W_{D^Q})\cap\mathcal{G}^{in}_n(A^{\tr})|$ & 0 & 0 & 0 & 2 & 286 \\
$|\mathcal{G}^{sp}_n(W_{A^{\tr}})\cap\mathcal{G}^{in}_n(A^{\tr})|$ & 0 & 0 & 0 & 2 & 286 \\
$|\mathcal{G}^{sp}_n(W_{D^{\deg}})\cap\mathcal{G}^{in}_n(A^{\tr})|$ & 0 & 0 & 0 & 2 & 286 \\
$|\mathcal{G}^{sp}_n(W_A)\cap\mathcal{G}^{in}_n(A^{\tr})|$ & 0 & 0 & 0 & 2 & 288 \\
$|\mathcal{G}^{sp}_n(W_Q)\cap\mathcal{G}^{in}_n(A^{\tr})|$ & 0 & 0 & 0 & 2 & 290 \\
$|\mathcal{G}^{sp}_n(A^{\tr})\cap\mathcal{G}^{sp}_n(W_D)|$ & 0 & 0 & 0 & 2 & 316 \\
\end{tabular}
    \caption{Number of graphs with up to $10$ vertices having a mate with respect to the spectrum or the SNF of two matrices associated with the graph at the same time. (part 2)}
    \label{tab:datos2}
\end{table}

\begin{table}[h]
    \footnotesize
    \centering
    \begin{tabular}{cccccc}
parameter & 6 & 7 & 8 & 9 & 10 \\
\hline
$|\mathcal{G}^{sp}_n(A)\cap\mathcal{G}^{sp}_n(W_{D^Q})|$ & 0 & 0 & 0 & 6 & 365 \\
$|\mathcal{G}^{sp}_n(A)\cap\mathcal{G}^{sp}_n(W_{A^{\tr}})|$ & 0 & 0 & 0 & 6 & 365 \\
$|\mathcal{G}^{sp}_n(A)\cap\mathcal{G}^{sp}_n(W_{D^{\deg}})|$ & 0 & 0 & 0 & 6 & 365 \\
$|\mathcal{G}^{sp}_n(A)\cap\mathcal{G}^{sp}_n(W_D)|$ & 0 & 0 & 0 & 6 & 373 \\
$|\mathcal{G}^{sp}_n(D^{\deg})\cap\mathcal{G}^{sp}_n(W_D)|$ & 0 & 0 & 0 & 6 & 373 \\
$|\mathcal{G}^{sp}_n(D^{\deg})\cap\mathcal{G}^{sp}_n(W_{D^Q})|$ & 0 & 0 & 0 & 6 & 375 \\
$|\mathcal{G}^{sp}_n(D^{\deg})\cap\mathcal{G}^{sp}_n(W_{A^{\tr}})|$ & 0 & 0 & 0 & 6 & 375 \\
\end{tabular}
    \caption{Number of graphs with up to $10$ vertices having a mate with respect to the spectrum or the SNF of two matrices associated with the graph at the same time. (part 3)}
    \label{tab:datos3}
\end{table}

\begin{table}[h]
    \footnotesize
    \centering
    \begin{tabular}{cccccc}
parameter & 6 & 7 & 8 & 9 & 10 \\
\hline
$|\mathcal{G}^{in}_n(A^{\tr})\cap\mathcal{G}^{in}_n(W_D)|$ & 0 & 2 & 12 & 147 & 6521 \\
$|\mathcal{G}^{sp}_n(A^{\tr})\cap\mathcal{G}^{in}_n(W_D)|$ & 0 & 2 & 36 & 151 & 6842 \\
$|\mathcal{G}^{in}_n(A^{\tr})\cap\mathcal{G}^{in}_n(W_A)|$ & 0 & 2 & 14 & 161 & 6764 \\
$|\mathcal{G}^{sp}_n(D^Q)\cap\mathcal{G}^{in}_n(W_D)|$ & 0 & 2 & 44 & 166 & 7423 \\
$|\mathcal{G}^{sp}_n(A^{\tr})\cap\mathcal{G}^{in}_n(W_A)|$ & 0 & 2 & 42 & 167 & 7099 \\
$|\mathcal{G}^{in}_n(D^Q)\cap\mathcal{G}^{in}_n(W_D)|$ & 0 & 2 & 16 & 170 & 7081 \\
$|\mathcal{G}^{sp}_n(Q)\cap\mathcal{G}^{in}_n(W_D)|$ & 0 & 2 & 48 & 170 & 7531 \\
\end{tabular}
    \caption{Number of graphs with up to $10$ vertices having a mate with respect to the spectrum or the SNF of two matrices associated with the graph at the same time. (part 4)}
    \label{tab:datos4}
\end{table}

\begin{table}[h]
    \footnotesize
    \centering
    \begin{tabular}{cccccc}
parameter & 6 & 7 & 8 & 9 & 10 \\
\hline
$|\mathcal{G}^{in}_n(L)\cap\mathcal{G}^{in}_n(W_{A^{\tr}})|$ & 0 & 2 & 31 & 630 & 20521 \\
$|\mathcal{G}^{in}_n(A^{\tr}_+)\cap\mathcal{G}^{in}_n(W_{A^{\tr}})|$ & 0 & 0 & 26 & 630 & 21410 \\
$|\mathcal{G}^{sp}_n(A)\cap\mathcal{G}^{in}_n(L)|$ & 0 & 6 & 32 & 637 & 24021 \\
$|\mathcal{G}^{sp}_n(D^{\deg}_+)\cap\mathcal{G}^{in}_n(W_{D^Q})|$ & 0 & 0 & 69 & 644 & 26026 \\
$|\mathcal{G}^{in}_n(L)\cap\mathcal{G}^{in}_n(W_{D^Q})|$ & 4 & 6 & 47 & 646 & 23219 \\
$|\mathcal{G}^{in}_n(D^L)\cap\mathcal{G}^{in}_n(W_{D^Q})|$ & 0 & 0 & 35 & 648 & 24478 \\
$|\mathcal{G}^{sp}_n(A^{\tr}_+)\cap\mathcal{G}^{in}_n(W_{D^Q})|$ & 0 & 0 & 71 & 656 & 27371 \\
\end{tabular}
    \caption{Number of graphs with up to $10$ vertices having a mate with respect to the spectrum or the SNF of two matrices associated with the graph at the same time. (part 5)}
    \label{tab:datos5}
\end{table}

\begin{table}[h]
    \footnotesize
    \centering
    \begin{tabular}{cccccc}
parameter & 6 & 7 & 8 & 9 & 10 \\
\hline
$|\mathcal{G}^{sp}_n(W_D)\cap\mathcal{G}^{in}_n(W_{D^{\deg}})|$ & 4 & 16 & 120 & 1275 & - \\
$|\mathcal{G}^{sp}_n(W_D)\cap\mathcal{G}^{in}_n(W_{D^{\deg}_+})|$ & 4 & 16 & 120 & 1275 & - \\
$|\mathcal{G}^{sp}_n(W_{D^Q})\cap\mathcal{G}^{in}_n(W_A)|$ & 4 & 16 & 120 & 1275 & - \\
$|\mathcal{G}^{sp}_n(W_{A^{\tr}})\cap\mathcal{G}^{in}_n(W_A)|$ & 4 & 16 & 120 & 1275 & - \\
$|\mathcal{G}^{sp}_n(W_{A^{\tr}})\cap\mathcal{G}^{in}_n(W_D)|$ & 4 & 16 & 120 & 1275 & - \\
$|\mathcal{G}^{sp}_n(W_{D^{\deg}})\cap\mathcal{G}^{in}_n(W_A)|$ & 4 & 16 & 120 & 1275 & - \\
$|\mathcal{G}^{sp}_n(W_{D^{\deg}})\cap\mathcal{G}^{in}_n(W_D)|$ & 4 & 16 & 120 & 1275 & - \\
\end{tabular}
    \caption{Number of graphs with up to $10$ vertices having a mate with respect to the spectrum or the SNF of two matrices associated with the graph at the same time. (part 6)}
    \label{tab:datos6}
\end{table}

\begin{table}[h]
    \footnotesize
    \centering
    \begin{tabular}{cccccc}
parameter & 6 & 7 & 8 & 9 & 10 \\
\hline
$|\mathcal{G}^{sp}_n(W_{D^{\deg}})\cap\mathcal{G}^{in}_n(W_{D^L})|$ & 4 & 16 & 120 & 1278 & - \\
$|\mathcal{G}^{sp}_n(W_{D^{\deg}})\cap\mathcal{G}^{in}_n(W_{A^{\tr}})|$ & 4 & 16 & 120 & 1278 & - \\
$|\mathcal{G}^{sp}_n(W_{D^{\deg}})\cap\mathcal{G}^{in}_n(W_{D^{\deg}})|$ & 4 & 16 & 120 & 1278 & - \\
$|\mathcal{G}^{sp}_n(W_D)\cap\mathcal{G}^{sp}_n(W_{D^Q})|$ & 4 & 16 & 124 & 1279 & 28577 \\
$|\mathcal{G}^{sp}_n(W_L)\cap\mathcal{G}^{sp}_n(W_D)|$ & 4 & 16 & 124 & 1279 & 28784 \\
$|\mathcal{G}^{sp}_n(W_D)\cap\mathcal{G}^{sp}_n(W_{D^L})|$ & 4 & 16 & 124 & 1279 & 28784 \\
$|\mathcal{G}^{sp}_n(W_A)\cap\mathcal{G}^{in}_n(W_{A^{\tr}_+})|$ & 4 & 16 & 125 & 1279 & 28919 \\
\end{tabular}
    \caption{Number of graphs with up to $10$ vertices having a mate with respect to the spectrum or the SNF of two matrices associated with the graph at the same time. (part 7)}
    \label{tab:datos7}
\end{table}

\begin{table}[h]
    \footnotesize
    \centering
    \begin{tabular}{cccccc}
parameter & 6 & 7 & 8 & 9 & 10 \\
\hline
$|\mathcal{G}^{sp}_n(D^L)\cap\mathcal{G}^{in}_n(A^{\tr})|$ & 0 & 0 & 46 & 1405 & 44141 \\
$|\mathcal{G}^{sp}_n(A^{\tr}_+)\cap\mathcal{G}^{in}_n(A^{\tr})|$ & 0 & 0 & 46 & 1405 & 44147 \\
$|\mathcal{G}^{sp}_n(D^Q)\cap\mathcal{G}^{in}_n(A^{\tr}_+)|$ & 0 & 0 & 46 & 1411 & 45348 \\
$|\mathcal{G}^{sp}_n(L)\cap\mathcal{G}^{in}_n(D^Q)|$ & 0 & 0 & 36 & 1421 & 44398 \\
$|\mathcal{G}^{sp}_n(Q)\cap\mathcal{G}^{in}_n(D^L)|$ & 0 & 0 & 42 & 1439 & 44848 \\
$|\mathcal{G}^{in}_n(D^L)\cap\mathcal{G}^{in}_n(D^Q)|$ & 0 & 0 & 44 & 1447 & 46239 \\
$|\mathcal{G}^{sp}_n(D^Q)\cap\mathcal{G}^{in}_n(D^L)|$ & 0 & 0 & 50 & 1459 & 46941 \\
\end{tabular}
    \caption{Number of graphs with up to $10$ vertices having a mate with respect to the spectrum or the SNF of two matrices associated with the graph at the same time. (part 8)}
    \label{tab:datos8}
\end{table}

\begin{table}[h]
    \footnotesize
    \centering
    \begin{tabular}{cccccc}
parameter & 6 & 7 & 8 & 9 & 10 \\
\hline
$|\mathcal{G}^{in}_n(D^{\deg}_+)\cap\mathcal{G}^{in}_n(W_{D^Q})|$ & 4 & 30 & 234 & 2057 & - \\
$|\mathcal{G}^{sp}_n(A^{\tr}_+)\cap\mathcal{G}^{in}_n(D)|$ & 0 & 0 & 62 & 2058 & 102255 \\
$|\mathcal{G}^{in}_n(D^{\deg}_+)\cap\mathcal{G}^{in}_n(W_{A^{\tr}})|$ & 4 & 30 & 224 & 2071 & - \\
$|\mathcal{G}^{in}_n(D^{\deg})\cap\mathcal{G}^{in}_n(D^{\deg}_+)|$ & 0 & 2 & 63 & 2118 & 57362 \\
$|\mathcal{G}^{in}_n(D^{\deg}_+)\cap\mathcal{G}^{in}_n(W_{D^{\deg}})|$ & 6 & 30 & 226 & 2153 & - \\
$|\mathcal{G}^{in}_n(D^Q)\cap\mathcal{G}^{in}_n(W_{D^{\deg}})|$ & 0 & 6 & 74 & 2159 & - \\
$|\mathcal{G}^{in}_n(A^{\tr})\cap\mathcal{G}^{in}_n(W_{D^Q})|$ & 0 & 6 & 66 & 2166 & - \\
\end{tabular}
    \caption{Number of graphs with up to $10$ vertices having a mate with respect to the spectrum or the SNF of two matrices associated with the graph at the same time. (part 9)}
    \label{tab:datos9}
\end{table}

\begin{table}[h]
    \footnotesize
    \centering
    \begin{tabular}{cccccc}
parameter & 6 & 7 & 8 & 9 & 10 \\
\hline
$|\mathcal{G}^{sp}_n(D^{\deg}_+)\cap\mathcal{G}^{in}_n(D)|$ & 0 & 18 & 208 & 4895 & 136058 \\
$|\mathcal{G}^{in}_n(D^{\deg}_+)\cap\mathcal{G}^{in}_n(W_D)|$ & 10 & 46 & 513 & 4981 & - \\
$|\mathcal{G}^{sp}_n(W_{A^{\tr}_+})\cap\mathcal{G}^{in}_n(A^{\tr})|$ & 4 & 20 & 188 & 5037 & - \\
$|\mathcal{G}^{sp}_n(W_{D^{\deg}_+})\cap\mathcal{G}^{in}_n(A^{\tr})|$ & 4 & 20 & 188 & 5037 & - \\
$|\mathcal{G}^{sp}_n(W_{A^{\tr}_+})\cap\mathcal{G}^{in}_n(D^Q)|$ & 4 & 20 & 194 & 5558 & - \\
$|\mathcal{G}^{sp}_n(W_{D^{\deg}_+})\cap\mathcal{G}^{in}_n(D^Q)|$ & 4 & 20 & 194 & 5558 & - \\
$|\mathcal{G}^{sp}_n(Q)\cap\mathcal{G}^{in}_n(W_Q)|$ & 0 & 28 & 295 & 5588 & - \\
\end{tabular}
    \caption{Number of graphs with up to $10$ vertices having a mate with respect to the spectrum or the SNF of two matrices associated with the graph at the same time. (part 10)}
    \label{tab:datos10}
\end{table}

\begin{table}[h]
    \footnotesize
    \centering
    \begin{tabular}{cccccc}
parameter & 6 & 7 & 8 & 9 & 10 \\
\hline
$|\mathcal{G}^{sp}_n(D^Q)\cap\mathcal{G}^{in}_n(W_{A^{\tr}_+})|$ & 6 & 38 & 414 & 6646 & - \\
$|\mathcal{G}^{sp}_n(D^Q)\cap\mathcal{G}^{in}_n(D^Q)|$ & 4 & 20 & 243 & 6676 & 255964 \\
$|\mathcal{G}^{sp}_n(Q)\cap\mathcal{G}^{in}_n(W_{A^{\tr}_+})|$ & 6 & 38 & 418 & 6676 & - \\
$|\mathcal{G}^{sp}_n(D^{\deg})\cap\mathcal{G}^{in}_n(W_{A^{\tr}_+})|$ & 6 & 38 & 395 & 6707 & - \\
$|\mathcal{G}^{sp}_n(Q)\cap\mathcal{G}^{in}_n(W_{D^{\deg}_+})|$ & 6 & 38 & 422 & 6728 & - \\
$|\mathcal{G}^{sp}_n(D^{\deg})\cap\mathcal{G}^{in}_n(W_{D^{\deg}_+})|$ & 6 & 38 & 399 & 6759 & - \\
$|\mathcal{G}^{sp}_n(Q)\cap\mathcal{G}^{in}_n(D^{\deg})|$ & 4 & 24 & 258 & 7018 & 262102 \\
\end{tabular}
    \caption{Number of graphs with up to $10$ vertices having a mate with respect to the spectrum or the SNF of two matrices associated with the graph at the same time. (part 11)}
    \label{tab:datos11}
\end{table}

\begin{table}[h]
    \footnotesize
    \centering
    \begin{tabular}{cccccc}
parameter & 6 & 7 & 8 & 9 & 10 \\
\hline
$|\mathcal{G}^{sp}_n(W_{A^{\tr}_+})\cap\mathcal{G}^{in}_n(D^L)|$ & 0 & 10 & 303 & 10982 & - \\
$|\mathcal{G}^{sp}_n(W_{D^{\deg}_+})\cap\mathcal{G}^{in}_n(D^L)|$ & 0 & 10 & 303 & 10982 & - \\
$|\mathcal{G}^{in}_n(D^{\deg})\cap\mathcal{G}^{in}_n(W_{A^{\tr}_+})|$ & 6 & 28 & 333 & 11507 & - \\
$|\mathcal{G}^{in}_n(D^{\deg})\cap\mathcal{G}^{in}_n(W_{D^{\deg}_+})|$ & 6 & 28 & 333 & 11561 & - \\
$|\mathcal{G}^{sp}_n(Q)\cap\mathcal{G}^{in}_n(Q)|$ & 8 & 59 & 603 & 13026 & 486524 \\
$|\mathcal{G}^{in}_n(A^{\tr}_+)\cap\mathcal{G}^{in}_n(W_{D^{\deg}_+})|$ & 0 & 14 & 382 & 13337 & - \\
$|\mathcal{G}^{in}_n(D^L)\cap\mathcal{G}^{in}_n(D^{\deg}_+)|$ & 0 & 12 & 332 & 13380 & 532587 \\
\end{tabular}
    \caption{Number of graphs with up to $10$ vertices having a mate with respect the spectrum or the SNF of two matrices associated with the graph at the same time. (part 12)}
    \label{tab:datos12}
\end{table}

\begin{table}[h]
    \footnotesize
    \centering
    \begin{tabular}{cccccc}
parameter & 6 & 7 & 8 & 9 & 10 \\
\hline
$|\mathcal{G}^{sp}_n(D)\cap\mathcal{G}^{in}_n(W_{D^{\deg}_+})|$ & 0 & 20 & 453 & 15822 & - \\
$|\mathcal{G}^{sp}_n(D^L)\cap\mathcal{G}^{in}_n(A^{\tr}_+)|$ & 0 & 16 & 456 & 15894 & 599752 \\
$|\mathcal{G}^{sp}_n(A^{\tr}_+)\cap\mathcal{G}^{in}_n(A^{\tr}_+)|$ & 0 & 16 & 456 & 15902 & 600515 \\
$|\mathcal{G}^{sp}_n(W_L)\cap\mathcal{G}^{in}_n(A^{\tr}_+)|$ & 0 & 16 & 456 & 15952 & - \\
$|\mathcal{G}^{sp}_n(W_{D^L})\cap\mathcal{G}^{in}_n(A^{\tr}_+)|$ & 0 & 16 & 456 & 15952 & - \\
$|\mathcal{G}^{in}_n(A^{\tr}_+)\cap\mathcal{G}^{in}_n(W_L)|$ & 0 & 16 & 456 & 15952 & - \\
$|\mathcal{G}^{in}_n(A^{\tr}_+)\cap\mathcal{G}^{in}_n(W_{D^L})|$ & 0 & 16 & 456 & 15952 & - \\
\end{tabular}
    \caption{Number of graphs with up to $10$ vertices having a mate with respect to the spectrum or the SNF of two matrices associated with the graph at the same time. (part 13)}
    \label{tab:datos13}
\end{table}

\begin{table}[h]
    \footnotesize
    \centering
    \begin{tabular}{cccccc}
parameter & 6 & 7 & 8 & 9 & 10 \\
\hline
$|\mathcal{G}^{sp}_n(A^{\tr}_+)\cap\mathcal{G}^{sp}_n(W_{D^L})|$ & 0 & 43 & 728 & 19757 & - \\
$|\mathcal{G}^{sp}_n(A^{\tr}_+)\cap\mathcal{G}^{in}_n(W_L)|$ & 0 & 43 & 728 & 19757 & - \\
$|\mathcal{G}^{sp}_n(A^{\tr}_+)\cap\mathcal{G}^{in}_n(W_{D^L})|$ & 0 & 43 & 728 & 19757 & - \\
$|\mathcal{G}^{sp}_n(D^L)\cap\mathcal{G}^{sp}_n(W_L)|$ & 0 & 43 & 745 & 20455 & - \\
$|\mathcal{G}^{sp}_n(D^L)\cap\mathcal{G}^{sp}_n(W_{D^L})|$ & 0 & 43 & 745 & 20455 & - \\
$|\mathcal{G}^{sp}_n(D^L)\cap\mathcal{G}^{in}_n(W_L)|$ & 0 & 43 & 745 & 20455 & - \\
$|\mathcal{G}^{sp}_n(D^L)\cap\mathcal{G}^{in}_n(W_{D^L})|$ & 0 & 43 & 745 & 20455 & - \\
\end{tabular}
    \caption{Number of graphs with up to $10$ vertices having a mate with respect the spectrum or the SNF of two matrices associated with the graph at the same time. (part 14)}
    \label{tab:datos14}
\end{table}

\begin{table}[h]
    \footnotesize
    \centering
    \begin{tabular}{cccccc}
parameter & 6 & 7 & 8 & 9 & 10 \\
\hline
$|\mathcal{G}^{sp}_n(W_{A^{\tr}_+})\cap\mathcal{G}^{in}_n(W_Q)|$ & 53 & 308 & 3249 & 55895 & - \\
$|\mathcal{G}^{sp}_n(W_{D^{\deg}_+})\cap\mathcal{G}^{in}_n(W_Q)|$ & 53 & 308 & 3249 & 55895 & - \\
$|\mathcal{G}^{in}_n(W_{D^Q})\cap\mathcal{G}^{in}_n(W_{D^{\deg}})|$ & 57 & 307 & 3296 & 56896 & - \\
$|\mathcal{G}^{in}_n(W_{D^Q})\cap\mathcal{G}^{in}_n(W_{A^{\tr}})|$ & 57 & 305 & 3286 & 56918 & - \\
$|\mathcal{G}^{in}_n(D)\cap\mathcal{G}^{in}_n(D^{\deg}_+)|$ & 6 & 112 & 1601 & 56995 & 3675977 \\
$|\mathcal{G}^{in}_n(W_Q)\cap\mathcal{G}^{in}_n(W_{A^{\tr}})|$ & 58 & 312 & 3319 & 57007 & - \\
$|\mathcal{G}^{in}_n(W_{D^Q})\cap\mathcal{G}^{in}_n(W_{A^{\tr}_+})|$ & 57 & 309 & 3376 & 57072 & - \\
\end{tabular}
    \caption{Number of graphs with up to $10$ vertices having a mate with respect to the spectrum or the SNF of two matrices associated with the graph at the same time. (part 15)}
    \label{tab:datos15}
\end{table}

\begin{table}[h]
    \footnotesize
    \centering
    \begin{tabular}{cccccc}
parameter & 6 & 7 & 8 & 9 & 10 \\
\hline
$|\mathcal{G}^{in}_n(D^{\deg}_+)\cap\mathcal{G}^{in}_n(W_{A^{\tr}_+})|$ & 46 & 397 & 4406 & 94142 & - \\
$|\mathcal{G}^{in}_n(D^{\deg}_+)\cap\mathcal{G}^{in}_n(W_{D^{\deg}_+})|$ & 46 & 397 & 4410 & 94817 & - \\
$|\mathcal{G}^{sp}_n(W_{D^{\deg}_+})\cap\mathcal{G}^{in}_n(A)|$ & 62 & 388 & 4584 & 95267 & - \\
$|\mathcal{G}^{sp}_n(W_{A^{\tr}_+})\cap\mathcal{G}^{in}_n(A)|$ & 62 & 388 & 4584 & 95277 & - \\
$|\mathcal{G}^{in}_n(A)\cap\mathcal{G}^{in}_n(W_D)|$ & 44 & 367 & 4774 & 96841 & - \\
$|\mathcal{G}^{in}_n(L)\cap\mathcal{G}^{in}_n(D)|$ & 56 & 428 & 4878 & 100594 & 4348954 \\
$|\mathcal{G}^{sp}_n(W_{A^{\tr}_+})\cap\mathcal{G}^{sp}_n(W_{D^{\deg}_+})|$ & 70 & 455 & 5390 & 106925 & - \\
\end{tabular}
    \caption{Number of graphs with up to $10$ vertices having a mate with respect to the spectrum or the SNF of two matrices associated with the graph at the same time. (part 16)}
    \label{tab:datos16}
\end{table}

\begin{table}[h]
    \footnotesize
    \centering
    \begin{tabular}{cccccc}
parameter & 6 & 7 & 8 & 9 & 10 \\
\hline
$|\mathcal{G}^{in}_n(Q)\cap\mathcal{G}^{in}_n(W_L)|$ & 78 & 620 & 7962 & 201282 & - \\
$|\mathcal{G}^{in}_n(Q)\cap\mathcal{G}^{in}_n(W_{D^L})|$ & 78 & 620 & 7962 & 201282 & - \\
$|\mathcal{G}^{in}_n(A)\cap\mathcal{G}^{in}_n(W_A)|$ & 66 & 576 & 8668 & 201402 & - \\
$|\mathcal{G}^{sp}_n(W_L)\cap\mathcal{G}^{in}_n(D^{\deg}_+)|$ & 46 & 495 & 7169 & 209822 & - \\
$|\mathcal{G}^{sp}_n(W_{D^L})\cap\mathcal{G}^{in}_n(D^{\deg}_+)|$ & 46 & 495 & 7169 & 209822 & - \\
$|\mathcal{G}^{in}_n(D^{\deg}_+)\cap\mathcal{G}^{in}_n(W_L)|$ & 46 & 495 & 7169 & 209822 & - \\
$|\mathcal{G}^{in}_n(D^{\deg}_+)\cap\mathcal{G}^{in}_n(W_{D^L})|$ & 46 & 495 & 7169 & 209822 & - \\
\end{tabular}
    \caption{Number of graphs with up to $10$ vertices having a mate with respect the spectrum or the SNF of two matrices associated with the graph at the same time. (part 17)}
    \label{tab:datos17}
\end{table}

\end{document}